\documentclass[12pt]{article}

\begin{document}

\begin{center}
{\LARGE\bf  Leibniz, Randomness \& the Halting Probability\vspace{1mm}}
\\ 
{\scriptsize\bf
Gregory Chaitin, IBM Watson Research Center, Yorktown Heights\vspace{2mm}}
\\ 
\emph{\textbf{Dedicated to Alan Turing on the 50th Anniversary of his Death}}
\end{center}

\bigskip

Turing's remarkable 1936 paper \emph{On computable numbers, with an application
to the Entscheidungsproblem}
marks a dramatic turning point in modern mathematics. On the one hand, the
computer enters center stage as a major mathematical concept.  On the other
hand, Turing establishes a link between mathematics and physics by talking about what
\emph{a machine} can accomplish.
It is amazing how far these ideas have come in a comparatively short amount of time;
a small stream has turned into a major river.
    
I have recently completed a small book about some of these developments, \emph{Meta Math!}
It is currently available as an e-book on my personal website, and also from arxiv.org,
and is scheduled to be published next year.  
Here I will merely give a few highlights.
     
My story begins with Leibniz in 1686, the year before Newton published his \emph{Principia.}
Due to a snow storm, Leibniz is forced to take a break in his attempts to improve the
water pumps for some important German silver mines, and writes down an outline of some
of his ideas, now known to us as the \emph{Discours de m\'etaphysique.}
Leibniz then sends a summary of the major points through a mutual friend to the 
famous fugitive
French \emph{philosophe} Arnauld, who is so horrified at what he reads that Leibniz never sends him nor 
anyone else the entire manuscript.
It languishes among Leibniz's voluminous personal papers and is only discovered and published 
many years
after Leibniz's death.
   
In sections V and VI of the \emph{Discours de m\'etaphysique,}
Leibniz discusses the crucial question of how we can distinguish a world in which science
applies from one in which it does not. Imagine, he says, that someone has splattered
a piece of paper with ink spots determining in this manner a finite set of points on the page. 
Nevertheless, Leibniz observes, there will always be a mathematical
equation that passes through this finite set of points.  Indeed, many good ways to
do this are now known. For example, what is called Lagrangian interpolation will do.
    
So the existence of a mathematical curve passing through a set of points cannot enable us to 
distinguish
between points that are chosen at random and those that obey some kind of a scientific law.
How then can we tell the difference? Well, says Leibniz,
if the equation must be extremely complex (``fort compos\'ee'') that is not a valid scientific law
and the points are random (``irr\'egulier'').  
     
Leibniz had a million other interests and earned a living as a consultant to princes, and as
far as I know after having this idea he never returned to this subject.  Indeed, he was
always tossing out good ideas, but rarely, with the notable exception of the infinitesimal
calculus, had the time to develop them in depth. 
    
The next person to take up this subject, as far as I know, is Hermann Weyl in his 1932 
book \emph{The Open World,} consisting of three lectures on metaphysics that Weyl gave at Yale
University.
In fact, I discovered Leibniz's work on complexity and randomness by reading this little book
by Weyl.
And Weyl points out that Leibniz's way of distinguishing between points that are random 
and those that follow a law by invoking the complexity of a mathematical formula is unfortunately
not too well defined, since it depends on what 
primitive functions 
you are allowed to use in writing that formula and
therefore varies as a function of time.
   
Well, the field that I invented in 1965 and which I call algorithmic information theory provides
a possible solution for the problem noticed by Hermann Weyl.
This theory defines a string of bits to be random, irreducible, structureless, if it is 
algorithmically incompressible, that is to say, if the size of the smallest computer
program that produces that particular finite string of bits as output is about the same
size as the output it produces.
    
So we have added two ideas to Leibniz's 1686 proposal.  First, we measure complexity in terms
of bits of information, i.e.,  0s and 1s.   Second, instead of mathematical equations,
we use binary computer programs.  
Crucially, this enables us to compare the complexity of a scientific theory (the computer program)
with the complexity of the data that it explains (the output of the computer program).
    
As Leibniz observed, for any data there is always a complicated theory, which is a computer program
that is the same size as the data. But that doesn't count.  It is only a real theory
if there is compression, if the program is much smaller than its output, both measured in 0/1 bits.
And if there can be no proper theory, then the bit string is algorithmically random or irreducible.
That's how you define a random string in algorithmic information theory.
    
I should point out that Leibniz had the two key ideas that you need to get this modern definition
of
randomness, he just never made the connection.  For Leibniz produced one of the first calculating
machines, which he displayed at the Royal Society in London, and he was also one of the first
people to appreciate base-two binary arithmetic and the fact that everything can be represented
using only 0s and 1s.  So,  as Martin Davis argues in his book \emph{The Universal Computer: The Road
from Leibniz
to Turing,} Leibniz was the first computer scientist, and he was also the first information 
theorist.  I am sure that Leibniz would have instantly understood and appreciated the modern definition of
randomness.
    
I should mention that A. N. Kolmogorov also proposed this definition of randomness.
He and I did this independently in 1965. 
Kolmogorov was at the end of his career, and I was a teenager at the
beginning of my own career as a mathematician.  
As far as I know, neither of us was aware of the Leibniz \emph{Discours.}
Let me compare Kolmogorov's work in this area with my own.  
I think that there are two key points to note.
    
Firstly,
Kolmogorov never realised as I did that our original definition of randomness was 
incorrect. It was a good initial idea but it was technically flawed. 
Nine years after he and I independently proposed this definition, I realised that
it was important for the computer programs that are used in the theory to be what I call
``self-delimiting''; without this it is not even possible to define my $\Omega$ number that I'll
discuss below. And there are other important changes that I had to make in the original definitions
that Kolmogorov never realised were necessary.
    
Secondly, Kolmogorov thought that the key application of these ideas was to be to obtain a new, algorithmic
version of probability theory.
It's true, that can be done, but it's not very interesting, it's too systematic a re-reading
of standard probability theory.  In fact, every statement that is true with probability one, merely
becomes a statement that must necessarily be true, for sure, for what are defined to be the random
infinite sequences of bits.
Kolmogorov never realised as I did that the really important application of these ideas was the new
light that they shed on G\"odel's incompleteness theorem and on Turing's halting problem.
    
So let me tell you about that now, and I'm sure that Turing would have loved these ideas if his
premature death had not prevented him from learning about them.
I'll tell you how my $\Omega$ number, which is defined to be the halting probability of a binary
program whose bits are generated using independent tosses of a fair coin, shows that in a sense
there is randomness in pure mathematics.
    
Instead of looking at individual instances of Turing's famous halting problem, you just 
put all possible computer programs into a bag, shake it well, pick out a program, and ask
what is the probability that it will eventually halt.  That's how you define the halting probability
$\Omega,$ and for this to work it's important that the programs have to be self-delimiting.
Otherwise the halting probability diverges to infinity instead of being a real number between
zero and one like all probabilities have to be. You'll have to take my word for this; I can't
explain this in detail here.
    
Anyway, once you do things properly you can define a halting probability $\Omega$ between zero and one.
The particular value of $\Omega$ that you get depends on your choice of computer programming language,
but its surprising properties don't depend on that choice.
    
And what is $\Omega$'s most surprising property? It's the fact that if you write $\Omega$ in binary, the
bits in its base-two expansion, the bits after the binary decimal point, seem to have absolutely
no mathematical structure. Even though $\Omega$ has a simple mathematical definition, its individual
bits seem completely patternless. In fact, they are maximally unknowable, they have, so to speak,
maximum entropy.  Even though they are precisely defined once you specify the programming language,
the individual bits are maximally unknowable, maximally irreducible. They seem to be mathematical facts that 
are true for no reason.
    
Why? Well, it is impossible to compress $N$ bits of $\Omega$ into a computer
program that is substantially smaller than $N$ bits in size (so that $\Omega$ satisfies the
definition of randomness of algorithmic information theory).  But not only does computation fail
to compress $\Omega,$
reason fails as well.
No formal mathematical theory whose axioms have less than $N$ bits of complexity can
enable us to determine $N$ bits of $\Omega.$  In other words, essentially the only way to be able to
prove what the values of $N$ bits of $\Omega$ are, is to assume what you want to prove as an axiom, which of
course is cheating and doesn't really count, because you are not using reasoning at all.
However, in the case of $\Omega,$ that is the best that you can ever do!
   
So this is an area in which mathematical truth has absolutely no structure, no structure that
we will ever be able to appreciate in detail, only statistically.  The best way of thinking about
the bits of $\Omega$ is to say that each bit has probability 1/2 of being zero and probability 1/2 of
being one, even though each bit is mathematically determined.
    
So that's where Turing's halting problem has led us, to the discovery of pure randomness
in a part of mathematics.
I think that Turing and Leibniz would be delighted at this remarkable turn of events.
   
Now I'd like to make a few comments about what I see as the philosophical
implications of all of this.
These are just my views, and they are quite controversial.
For example, even though
a recent critical review of two of my books in the \emph{Notices of the American Mathematical Society}
does not claim that there are any technical mistakes in my work, the reviewer strongly
disagrees with my philosophical conclusions, and in fact he claims that my work has no
philosophical implications whatsoever.
So these are just my views, they are certainly not a community  consensus, not at all.
    
My view is that $\Omega$ is a much more disagreeable instance of mathematical incompleteness
than the one found by G\"odel in 1931, and that it therefore forces our hand philosophically.
In what way? Well, in my opinion, in a quasi-empirical direction, which is a phrase coined
by Imre Lakatos when he was doing philosophy in England after leaving Hungary in 1956.
In my opinion, $\Omega$ suggests that even though math and physics are different, perhaps they
are not as different as most people think.
    
What do I mean by this?
(And whether Lakatos would agree or not, I cannot say.)
I think that physics enables us to compress our experimental data, and math enables us to
compress the results of our computations, 
into scientific or mathematical theories as the case may be.  
And I think that neither math nor science gives absolute certainty; that is an asymptotic
limit unobtainable by mortal beings.
And in this connection I should mention the book (actually two books) just published by
Borwein and Bailey on experimental math.  
    
To put it bluntly, if the incompleteness phenomenon
discovered by G\"odel in 1931 is really serious---and I believe that Turing's work and my
own work suggest that incompleteness is much more serious than people think---then perhaps mathematics
should be pursued somewhat more in the spirit of experimental science rather than always 
demanding proofs for everything. In fact, that is what theoretical computer scientists
are currently doing.  Although they may not want to admit it, and refer to 
\emph{\textbf{P}} $\neq$ \emph{\textbf{NP}}
as an unproved hypothesis, that community is in fact behaving as if this were a new axiom,
the way that physicists would.
   
At any rate, that's the way things seem to me. Perhaps by the time we reach the centenary
of Turing's death this quasi-empirical view will have made some headway, 
or perhaps instead these foreign ideas will be utterly
rejected by the immune system of the math community. For now they certainly are rejected.
But the past fifty years have brought us many surprises, 
and I expect that the next fifty years will too,
a great many indeed.
     
\section*{References}
     
\begin{enumerate}
\item
Chapter on Leibniz in E. T. Bell, \emph{Men of Mathematics,} Simon \& Schuster, 1986.
\item
R. C. Sleigh, Jr., \emph{Leibniz and Arnauld: A Commentary on Their Correspondence,} 
Yale University Press, 1990.
\item
H. Weyl, \emph{The Open World: Three Lectures on the Metaphysical Implications of Science,} 
Yale University Press, 1932,
Ox Bow Press, 1994. 
\item
Article by Lakatos in T. Tymoczko, \emph{New Directions in the Philosophy of Mathematics,}
Princeton University Press, 1998.
\item
M. Davis, \emph{The Universal Computer: The Road from Leibniz to Turing,} Norton, 2000.
\item
Play on Newton vs.\ Leibniz in C. Djerassi, D. Pinner, \emph{Newton's Darkness: Two Dramatic Views,}
Imperial College Press, 2003.
\item
J. Borwein, D. Bailey, \emph{Mathematics by Experiment: Plausible Reasoning in the 21st Century,}
A. K. Peters, 2004.
\\  
J. Borwein, D. Bailey, R. Girgensohn, \emph{Experimentation in Mathematics: Computational Paths
to Discovery,} A. K. Peters, 2004.
\item
G. Chaitin, \emph{Meta Math! The Quest for Omega,} 
\\  
http://cs.umaine.edu/\verb!~!chaitin/omega.html,
\\
http://arxiv.org/abs/math.HO/0404335. 
\\    
To be published by Pantheon Books in 2005.
\end{enumerate}

\end{document}